\documentclass[]{interact}
\usepackage{epstopdf}
\usepackage{subfigure}
\usepackage{url}

\usepackage[numbers,sort&compress]{natbib}
\bibpunct[, ]{[}{]}{,}{n}{,}{,}
\renewcommand\bibfont{\fontsize{10}{12}\selectfont}
\makeatletter
\def\NAT@def@citea{\def\@citea{\NAT@separator}}
\makeatother

\usepackage{xcolor} 
\usepackage{ulem}

\newcommand{\ansin}[1]{\textcolor{black}{#1}}

\usepackage{amsmath}
\DeclareMathOperator{\Tr}{Tr}

\theoremstyle{plain}

\theoremstyle{definition}

\theoremstyle{remark}

\newcommand{\dd}{\,{\rm d}} 

\begin{document}


\title{Optimal control of thermal and mechanical loads in activation processes of mechanical components}

\author{
\name{Nicolai Friedlich\textsuperscript{a}\thanks{CONTACT Nicolai Friedlich. Email: nicolai.friedlich@outlook.de}, Hanno Gottschalk\textsuperscript{b}, Georg Vossen\textsuperscript{a}}
\affil{
\textsuperscript{a}Institute of Modelling and High-Performance Computing, Niederrhein University of Applied Sciences, Krefeld, Germany; \\
\textsuperscript{b} Institute of Mathematics, Technical University Berlin, Germany}}

\maketitle

\begin{abstract}
This paper develops a mathematical framework that aims to control the temperature and rotational speed in the activation process of a gas turbine in an optimal way. These controls influence the deformation and the stress in the component due to centripetal loads and transient thermal stress. An  optimal control problem is formulated as the minimization of maximal \ansin{von }Mises stress over a given time and over the whole  component. To find a solution for this, we need to solve the linear thermoelasticity and the heat equations using the finite element method. The results for the optimal control as functions of the rotation speed and external gas temperature over time  are  computed by  sequential quadratic programming, where gradients are computed using finite differences.\ansin{The overall outcome reveals a significant reduction of approximately 10\%, from $830 \frac{N}{mm^2}$ to $750 \frac{N}{mm^2}$, in \ansin{von }Mises stress by controlling two parameters, along with the temporal separation of physical control phenomena.} \\
\end{abstract}


\section{Introduction}\label{sec1}
When technical facilities are activated or deactivated, the built-in components undergo load cycles that cause fatigue and, ultimately, failure. The safe operation of such facilities therefore requires service concepts to counteract the material fatigue by repair and replacement, but also concepts for the gentle operation of components that minimize the loads occurring. While this is general wisdom, so far little mathematical or computational research effort has been invested to minimize negative effects of activation cycles by controlling them actively and in an optimal way. In this article, we present a first study towards the optimal control of activation and deactivation processes. 

A technical device, where fatigue is highly relevant, is the gas turbine, which is mostly used in jet engines for aviation, but also for energy production. In both fields, activation and deactivation happens frequently, e.g.\ during start and landing, or when gas turbines are used in the energy grid  to balance volatile energy influx by renewable sources of energy production. For gas turbines, the loss of single turbine disks is contained in the turbine casing, but a rotor burst leads to a disintegration of the entire device with potentially catastrophic consequences. Therefore, the mechanical integrity of gas turbine rotor disks is crucial. At the same time, peak stresses in a gas turbine rotor do not necessarily occur under full load, but during the activation process itself. This is so, as the rotor usually is the heaviest part of the turbine with the largest thermal inertia. So, during activation, the exterior of the rotor is in contact with the hot working fluid and undergoes thermal expansion, whereas the interior part is still cool, which leads to transient thermal stress. At the same time, rotation speed and burner temperature are controllable for gas turbines, where constraints have to be respected for the duration of the activation process. This complex behavior makes gas turbines predestined for our study and consequently we work with an example motivated by the activation of a stationary gas turbine.     

The coupling of heat transfer and elastic deformation is a widely used phenomenon in many engineering applications where body parts are in varying temperature fields and have to withstand to mechanical loads. In \citep{Zas93}  an optimal control of the heating of a thermoelastic plate by internal heat sources is given. In \citep{Saa21} a thermoelastic sensitivity analysis has been made. From this work, the underlying equations were  adapted. In \citep{Zha16}, an online optimal control schemes of inlet steam temperature during start up of steam turbines  has been made. Contrary to this work, here we apply an optimal control with two control variables, the temperature and the rotation.
In \citep{ban17}, an optimization of a rotor with respect to the turbine load was carried out. 
Also turbine topology optimization was conducted in \citep{Jin02},\citep{Alk21} and \citep{Wan21} regarding to thermal and mechanical loads.
\ansin{In \citep{KHA2023_1}, \citep{KHA2023_2} and \citep{RAZ2023}, a detailed analysis for the heat transport rate was carried out. Here, we use a simplified boundary condition for the working fluid with a fixed convective heat transfer coefficient.}

In this paper, we will compute a solution of an optimal control problem for an activation process of a turbine regarding to thermal and mechanical loads. The mechanical as well as the thermal load will be the control variables over a given time and the goal is to reduce the maximum \ansin{von }Mises stress over the whole time. In general optimal control problem can be solved in two ways: First optimize then discretize (FOTD) or first discretize then optimize (FDTO) \citep{Tro15}. We will use the latter one. This means we will first discretize the volume with a 3D mesh, the time with implicit Euler and then solve the resulting finite-dimensional optimization problem.

The paper is organized as follows: Section 2 provides a detailed technical background and introduces the model formulation. Section 3 delves into the optimal control problem, discussing objectives, and constraints. In Section 4, a Finite Element Method (FEM) Simulation model is presented, comprising a weak formulation and an overview of preprocessing and implementation techniques. Section 5 explores the solution of the optimal control problem using Sequential Least Squares Programming (SLSQP) algorithm. Finally, in Section 6, a summary of the paper's findings is provided, along with an outlook on potential future research.

\section{Technical background and model formulation}\label{sec2}
Gas turbines play an essential role as back up capacities in energy grids with volatile influx of renewable energy from solar panels or wind turbines. This makes it necessary to start gas turbines quickly in order to meet the consumer's demands. 

During the starting procedure, the gas turbine burner heats up the flowpath, while the gas turbine starts turning. While the rotation leads to centripetal mechanical loads, the non uniform distribution of heat during the activation process results in thermal loads. Both mechanisms jointly cause mechanical fatigue which effectively limits the life of the gas turbine. For safe operation it is therefore desirable, to minimize peak stress in important components via a controlled activation procedure that increases the angular velocity and the temperature in an optimized manner.

Thermal inertia during heating leads to the aforementioned non uniform temperature distribution and thereby, potentially, to transient peak stress. This mechanism is more severe for massive parts, like the gas turbine rotor, and is less important for parts with little mass, as gas turbine blades and vanes. At the same time, the gas turbine rotor disks are the most critical parts of the turbine, as a rotor burst is not contained by the casing. We therefore choose a simplified model of a rotor disk as the most relevant application for our optimized activation half-cycle.

\begin{figure}[t]
\centering
\subfigure[Full 3D model of a simplified turbine. Outer diameter is about 1.4 m and represents a larger stage of a turbine. This stage is especially interesting for thermal and rotational loading.]{%
\resizebox*{6cm}{!}{\includegraphics{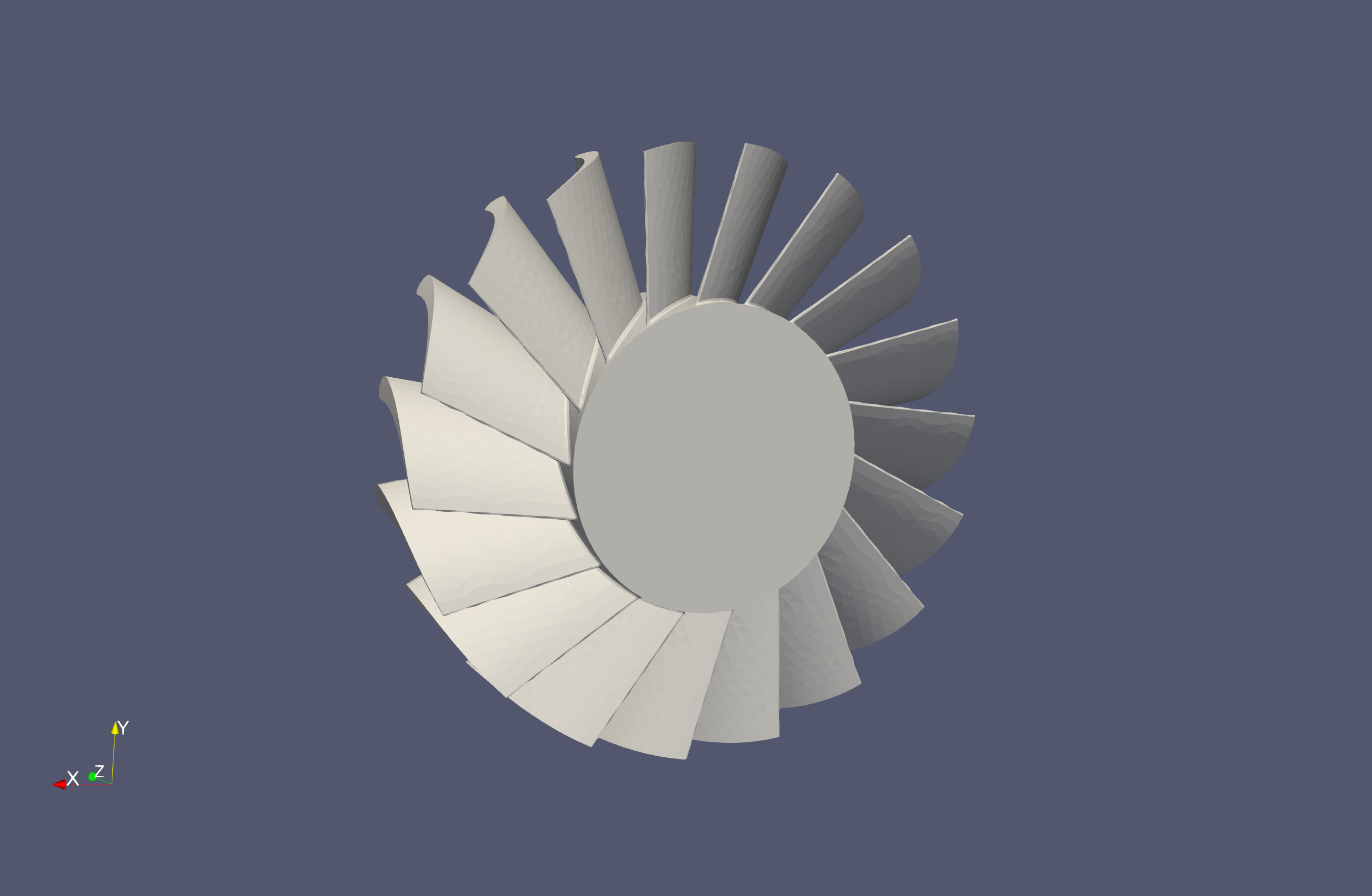}}}\hspace{5pt}
\subfigure[Due to axis symmetry reduced model with only one turbine blade. The maximal load is expected to be in the transition from turbine blade to rotor. This is the calculation volume $\Omega$.]{
\resizebox*{6cm}{!}{\includegraphics{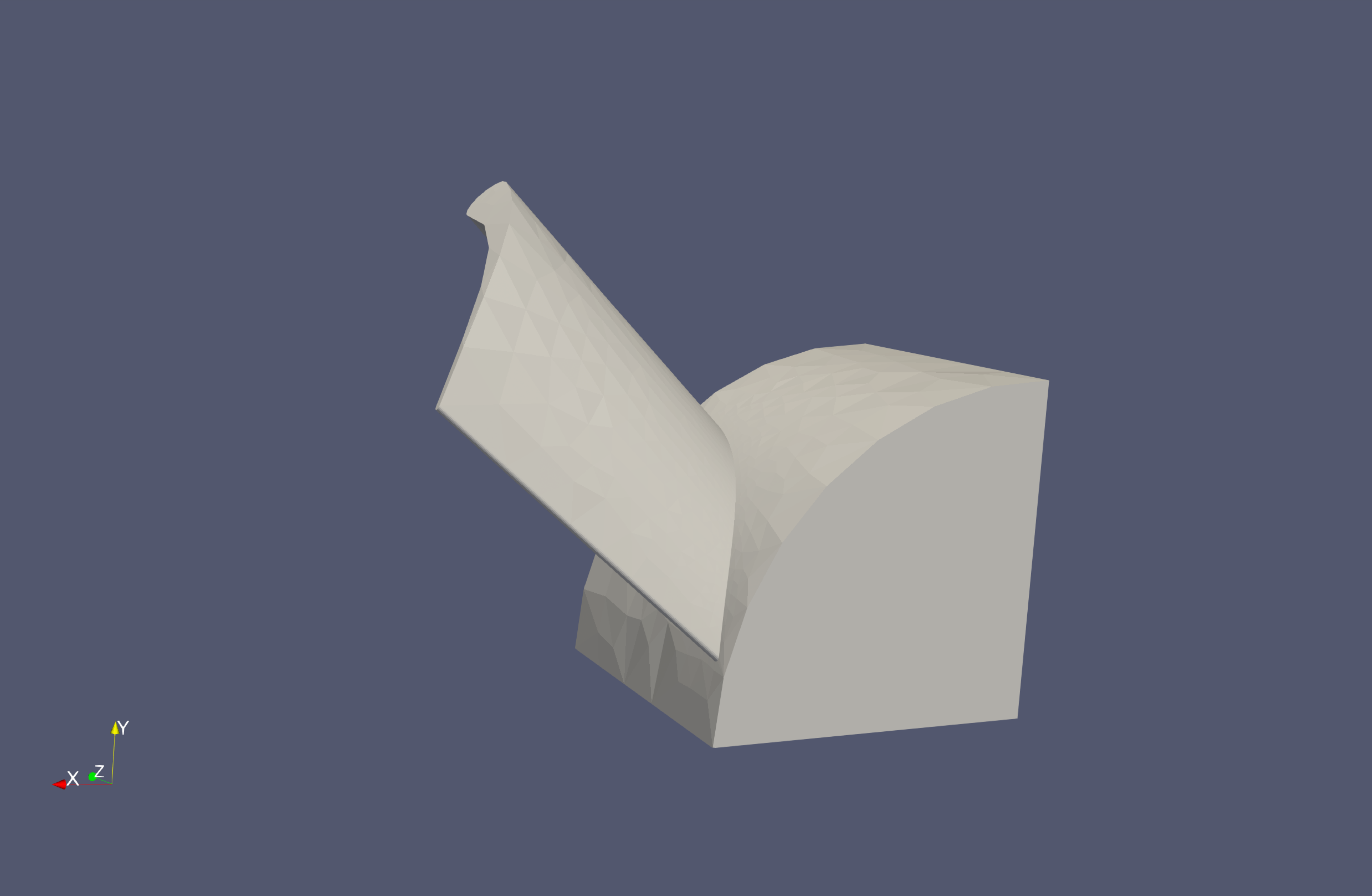}}}
\caption{Full (left, (a)) and reduced (right, (b)) 3D  model.} 
\label{fig:3d-figure}
\end{figure}
Let us now present the details of our mathematical model. \ansin{Throughout the paper, matrices are represented by bold uppercase letters, vectors by bold lowercase letters, and scalars by lowercase normal letters}. Let $\Omega\subseteq \mathbb{R}^3$ be the portion of the three dimensional space that is filled by the component, cf.\ Figure \ref{fig:3d-figure} and let $\tau = [0,t_f] $ be the underlying time sequence with a final time $t_f$. \ansin{The region $\Omega$ is tied to the Cartesian coordinate system and the front center of the cylinder in Figure \ref{fig:3d-figure} (b) is coordinate origin $(0, 0, 0)^T$}.   When an external load $\boldsymbol{f}:\tau\times\Omega\to\mathbb{R}^3$ is applied, the elastic material filling $\Omega$ undergoes a deformation that can be modeled by a vector field $\boldsymbol{u}=(u_x,u_y,u_z)^T:\tau\times\Omega\to\mathbb{R}^3$\ with $u_x$ denoting the deformation in $x$-direction, and $u_y$ and $u_z$ representing the deformations in $y$- and $z$-direction, respectively. Under the assumption of linear elasticity and within the context of linear small deformation theory, $\boldsymbol{u}$ fulfills the partial differential equation of thermoelasticity
\begin{align}\label{eq:elast}
&\nabla \cdot \boldsymbol{\Sigma}(\boldsymbol{u}) + \boldsymbol{f} = \boldsymbol{0} ~~~~~~~~~~~~~~~~~~~~~~~~~~~~~~~~~~~~~~~~~~~~~~~~~~~~ \text{in } \Omega
\end{align}
where
\begin{align}\label{eq:elast2}
&{\boldsymbol{\Sigma}}(\boldsymbol{u}) = 2 \mu\boldsymbol{\mathcal{E}}(\boldsymbol{u})+[\lambda \Tr(\boldsymbol{\mathcal{E}}(\boldsymbol{u})) -\alpha(3\lambda+2\mu)(T-T_0)]\boldsymbol{I} ~~~~~ \text{in } \Omega
\end{align}
is the stress tensor with $ \lambda,\mu>0$ the Lam\'e constants, $\boldsymbol{I}$ the $3\times3$ identity matrix and $T:\tau\times\Omega\to\mathbb{R}$ some given field of temperatures. $T_0$ is a reference temperature that corresponds to the temperature state of the non activated component. $\Tr$ denotes the trace operator. Thermal stresses are driven by $\alpha$,  thermal expansion coefficient, and the difference between $T(x)$, the temperature at $x\in\Omega$ and $T_0$. The elastic strain tensor is given by
\begin{align}\label{eq:strain}
\boldsymbol{\mathcal{E}}(\boldsymbol{u}) = \frac{1}{2}(\nabla \boldsymbol{u}+\nabla \boldsymbol{u}^\top)
\end{align}
and the loading due to centripetal acceleration due to rotation around the $z$ - axis of the axial turbo machine 
\begin{align}\label{eq:load}
\boldsymbol{f} = \rho \omega^2 \boldsymbol{x_p}, \quad \boldsymbol{x_p}=(x,y,0)^\top.
\end{align}
   The density is denoted as $\rho$, $\omega$ is the rotation and $\boldsymbol{x_p}$ are the projected radial components of the coordinate vector $(x,y,z)^\top\in\Omega$ where these rotation force apply to.

\begin{figure}[t]
\centerline{\includegraphics[width=0.7\linewidth]{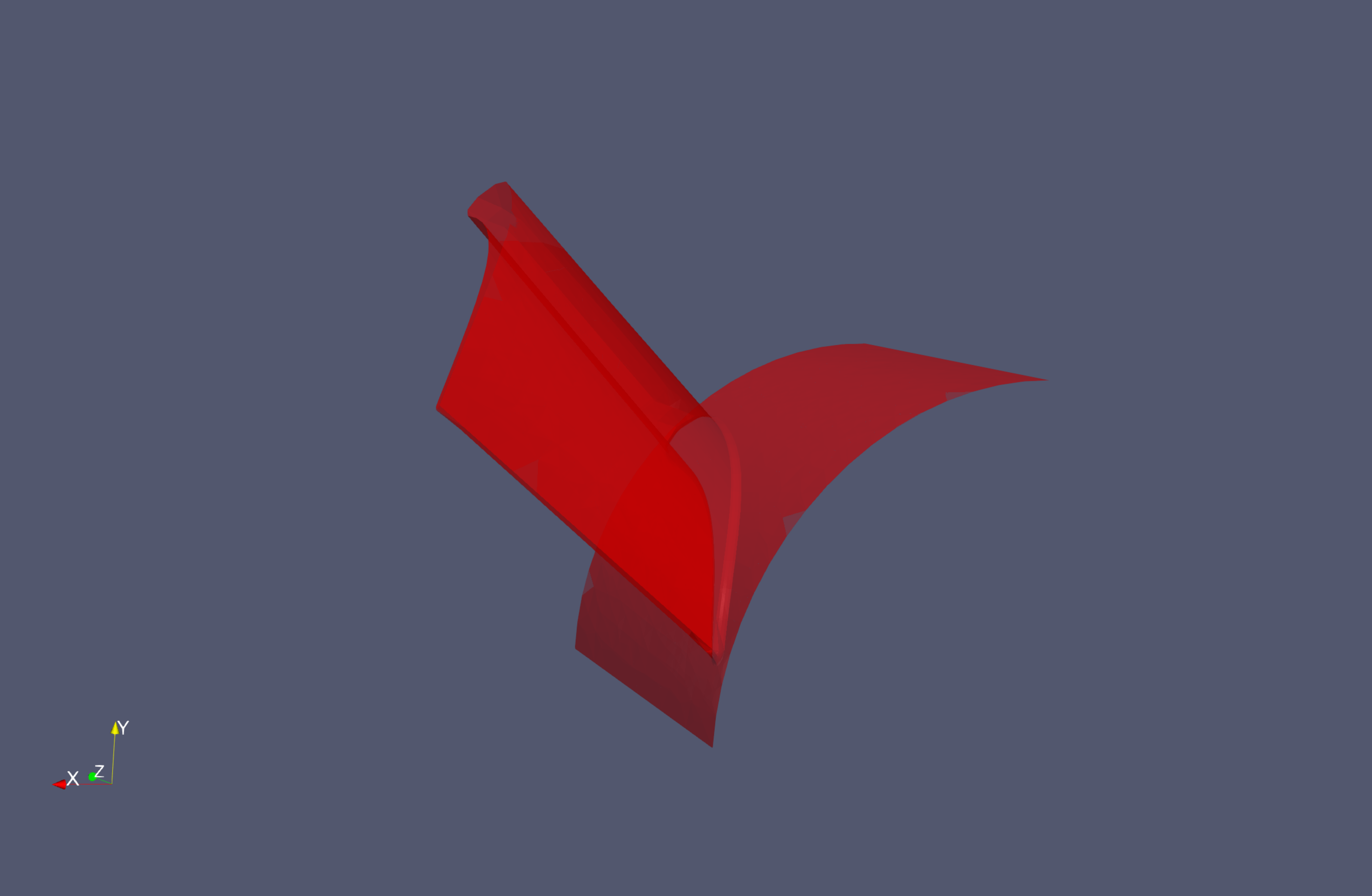}}
\caption{Boundary surface $(\partial\Omega)_R$ for external heating. The "R" subscript abbreviates for Robin boundary.}\label{fig:heat_bound}
\end{figure}

The temperature distribution $T :\tau\times \Omega\to\mathbb{R}$ is obtained as the solution of the time dependent heat equation
\begin{align}
\label{eq:heat}
\begin{split}
\rho c_p \frac{\partial T}{\partial t}-k\Delta T&= 0 ~\quad \text{in } \Omega;\\
T(t=0)&=T_0 \quad \text{in } \Omega;
\end{split}
\end{align}
where $\Delta=\frac{\partial^2}{\partial x^2}+\frac{\partial^2}{\partial y^2}+\frac{\partial^2}{\partial z^2}$ is the Laplace operator, $c_p$ the specific heat capacity and $k$ is the thermal conductivity. The boundary conditions at the surface where heat exchance with the heated working fluid is possible are given by Robin or convective boundary conditions
\begin{align}
\label{eq:heat_bound}
-k\frac{\partial T}{\partial \boldsymbol{n}} = h(T-T_\mathrm{e}) \quad \text{in } (\partial\Omega)_R
\end{align}
where $\frac{\partial }{\partial \boldsymbol{n}}$ stands for the directional derivative in the direction of the outward normal vector $\boldsymbol{n}$. The heat transfer coefficient $h>0$, is assumed to be fixed\footnote{A refined analysis should model $h$ as a function of the fluid boundary layer, which depends on temperature and mass flow of the working fluid.}. The notation $(\partial \Omega)_R$ denotes a Robin boundary and the surfaces, where this type of boundary conditions is applied, are depicted in Figure \ref{fig:heat_bound}. 

Note that the speed of rotation $\omega=\omega(t)$ and the external temperature $T_\mathrm{e}(t)$ are both variables that can be controlled during the activation process. This, together with the time dependent heat equation \eqref{eq:heat} induce time dependent solutions $ \boldsymbol{u}(t)$ of \eqref{eq:elast}.

To conclude this section, let us complement the boundary conditions for a computational domain that makes use of the component's rotation symmetry. In order to avoid the redundant numerical computation of temperature and stress in identical segments, we simplify the blade count to one and introduce plane symmetry boundary condition for the cylinder part on the 'internal' boundary sections of one specific segment, see Figure \ref{fig:3d-figure} (b), i.e.\ we set
\begin{align}
u_y(t,0,y,z)=0,\quad  u_x(t,x,0,z)=0, \label{eq:heat_bound2}\\
\frac{\partial T}{\partial \boldsymbol{n}}(t,0,y,z)=0, \quad \frac{\partial T}{\partial \boldsymbol{n}}(t,x,0,z)=0,\label{eq:heat_bound3}\\
\frac{\partial T}{\partial \boldsymbol{n}}(t,x,y,0)=0, \quad \frac{\partial T}{\partial \boldsymbol{n}}(t,x,y,z_l)=0.\label{eq:heat_bound4}
\end{align}
The boundary conditions \eqref{eq:heat_bound2} for $\boldsymbol{u}$ and \eqref{eq:heat_bound3}, \eqref{eq:heat_bound4} for $T$ are symmetric and null Neumann boundary conditions, where $z_l$ is the length of the part in $z$-direction. This means the deformation perpendicular to the $x$ and $y$ axis is zero as well as the part is isolated at any boundary except for the Robin boundary \eqref{eq:heat_bound}.
\ansin{The first boundary conditions in \eqref{eq:heat_bound2} and in \eqref{eq:heat_bound3} are applied along the bottom surface in Figure \ref{fig:3d-figure} (b). The second boundary conditions in \eqref{eq:heat_bound2} and in \eqref{eq:heat_bound3} are applied along the right hand side surface of the cut cylinder in Figure \ref{fig:3d-figure} (b). The boundary conditions in \eqref{eq:heat_bound4} are applied along the front and the back surface of the cylinder, respectively.}
We also do not apply any forces on the surface $\partial\Omega_D$ in Figure \ref{fig:heat_bound}.  In a later, refined modeling approach, pressure boundary conditions from the static pressure of the working fluid should be applied as well. However, this involves complex 3D fluid dynamics simulation and is beyond the scope of the present article. 
\section{Optimal control problem}\label{sec3}
The durability of turbine (as well as other engineering system) components is known to be enhanced if they are exposed to minimal stresses inside. Stress peaks can lead to unwanted deformations, porosities or even cracks in the parts. The strain of turbine components is particularly high during the starting process. 

We will therefore formulate an optimal control problem for the turbine activation. Let us consider the rotation  $\omega:\tau\to\mathbb{R}$ in \eqref{eq:load} and the external temperature  $T_e:\tau\to\mathbb{R}$ in \eqref {eq:heat_bound} as control functions in the fixed and given time domain $\tau=[0,t_f]$ of the starting process. The goal of the optimal control problem is to minimize maximal stresses considered over time and space, i.e., we minimize the cost functional
\begin{align} \label{eq:opti-cost}
J(T_\mathrm{e},\omega):=   \max_{t}(\max_{x,y,z}(\sigma_v)).
\end{align}
Here, $\sigma_v=\sigma_v(t,x,y,z)$ is the \ansin{von }Mises stress given by
\begin{align}\label{eq:Mises-stress}
\sigma_v &= \sqrt[]{\frac{3}{2} \boldsymbol{\Sigma^{\textrm{dev}}}:\boldsymbol{\Sigma^{\textrm{dev}}}}
\end{align}
where 
\begin{align}
\boldsymbol{\Sigma^{\textrm{dev}}} &= \boldsymbol{\Sigma} - \frac{1}{3} \textrm{Tr}(\boldsymbol{\Sigma}) \boldsymbol I
\end{align}
is the stress deviator tensor with $\boldsymbol\Sigma$ as defined in \eqref{eq:elast2}  \ansin{and} the matrix scalar product notation 
\begin{align*}
\boldsymbol{A}:\boldsymbol{B}=\textrm{Tr}(\boldsymbol{A}^\top \boldsymbol{B})=\sum_{i,j=1}^na_{ij}b_{ij}
\end{align*}
for two $n\times n$-matrices $\boldsymbol{A}$, $\boldsymbol{B}$ with entries $a_{ij}$, $b_{ij}$, $1\le i,j\le n$, has been used. From a technical point of view, we assume the control functions to be continuous. As discussed in the next paragraph, this will also be mathematically reasonable and furthermore, $\omega$  is assumed to be differentiable.

The starting process is determined by the three conditions that the external temperature $T_e$ has reached at least some desired temperature $T_{e,f}$, that the rotation $\omega$ has reached a desired rotation $\omega_f$ at the final time $t_f$ and that a the maximum temperature reached at least some desired condition. We hence postulate
\begin{align} \label{eq:opti-constr1}
T_e(t_{f})\geq T_{e,f}, \quad \omega(t_f)=\omega_f, \quad \max_{x,y,z}T(t_f)\geq T_f.
\end{align}
We can see that continuity of the control functions is essential in this formulation to prohibit solutions which jump from zero to the desired state in the final time $t_f$. Furthermore, we are interested in solutions with bounded rotation acceleration to prevent peaks in the rotation velocity, i.e., we set 
\begin{align}\label{eq:opti-constr2}
\frac{\partial \omega}{\partial t} \leq \omega_{lim},\quad t\in\tau
\end{align}
as pointwise constraint for the derivative of the control $\omega$. This assumption is also natural, as the rotational inertia forces will prevent an arbitrary acceleration of the angular velocity. 
In addition the control variable underlies some bounds
\begin{align}
\omega_{min}<\omega<\omega_{max},\label{eq:Te_bound}\\
 T_{e,min}<T_e< T_{e,max}.\label{eq:w_bound}
\end{align}
This bounds will start at $\omega_{min} = 0$ and  $T_{e,min} = 0$ and will go to a maximum of $\omega_{max} = 60 $ hz and $T_{e,max} = 1000 $ °C.
Finally, the dynamics of our model with boundary and initial conditions, as given in Section \ref{sec2}, are constraints of the problem.

To summarize, the optimal control problem reads:
\begin{align} \label{eq:opti}
\begin{split}
&\mbox{Minimize}\quad \eqref{eq:opti-cost}\\
&\text{subject to} \quad \eqref{eq:elast},\;\eqref{eq:heat},\;\eqref{eq:heat_bound2},\;\eqref{eq:heat_bound3},\;\eqref{eq:heat_bound4},\;\eqref{eq:opti-constr1},\;\eqref{eq:opti-constr2}\;\eqref{eq:Te_bound}\;\eqref{eq:w_bound}.
\end{split}
\end{align}

\section{FEM Simulation model}\label{sec4}
FEM is a method for solving partial differential equations and it is very widely used in the engineering field to solve structural mechanics problems. We are looking on a structural mechanics problem coupled with the heat equation. One of the main strengths is the flexibility to use complex computational domains, so we can calculate 3d models easily.
This section captures the used methods of this work. Firstly we will show the weak formulation of the governing equation which we need for the implementation afterwards. Then we show our 3D model and the FEM Model where we used FEniCS \citep{Log12} and \citep{Aln15}, an open source FEM solver. We also provide an overview about the workflow in the appendix because it is not trivial when getting started with complex 3d domains in FEM.

\subsection{Weak formulation}\label{subsec42}
To use FEM, we have to formulate the weak form of the governing equations. We adapted the formulation from \citep{Saa21} \ansin{and \citep{Hof85}} and added a transient part.
Applying the weak formulation on the heat equation \eqref{eq:heat} yields
\begin{align}
\int_{\Omega}\rho c_p \frac{\partial T}{\partial t} \cdot v_T \ \dd x-\int_{\Omega}k\Delta T \cdot v_T \ \dd x= 0
\end{align}
where $v_T \in H^1(\Omega)$ is a \ansin{test function} for the temperature and $\dd x$ denotes the integration over the whole domain $\Omega $. Hereby, $H^1$ is the Sobolev space of functions whose derivative is $L^2$ integrable; see, e.g., \citep{Tro15} for details. Using partial integration for the second term results in
\begin{align}
\int_{\Omega}\rho c_p \frac{\partial T}{\partial t} \cdot v_T \ \dd x+\int_{\Omega}k\nabla T \cdot \nabla v_T \ \dd x- k\int_{\partial \Omega} (k \nabla T \cdot \boldsymbol{n} )v_T \ \dd A = 0.
\end{align}
Here $dA$ denoted the integration over the a boundary $\partial\Omega$. With help of equation \eqref{eq:heat_bound}, the last term can be written  as
\begin{align}
- \int_{\partial \Omega} (k \nabla T \cdot \boldsymbol{n} )v_T \ \dd A =
- \int_{\partial \Omega} k\frac{\partial T}{\partial \boldsymbol{n}}v_T \ \dd A = 
- \int_{\partial \Omega} h(T-T_\mathrm{e})v_T \ \dd A
\end{align} 

The time derivatives are now replaced by an backward Euler scheme, so that the previous weak form at the time increment $n+1$ is now:
\begin{align}
\int_{\Omega}\rho c_p \dfrac{T_{n+1}-T_n}{\Delta t} \cdot v_T \ \dd x+\int_{\Omega}k\nabla T_{n+1} \cdot \nabla v_T \ \dd x= -\int_{\partial \Omega} h(T_{n+1}-T_\mathrm{e}))v_T \ \dd A
\end{align}
The same has to be done for the thermoelasticity equation \eqref{eq:elast}:
\begin{equation}
\int_{\Omega} \boldsymbol{\Sigma}(u):\boldsymbol{\mathcal{E}}(v_M) \ \dd x = \int_{\Omega} f \cdot v_M \ \dd x +\int_{\partial\Omega}v_M \cdot \boldsymbol{\Sigma}(u)\cdot n \ \dd A
\end{equation}
where $v_M \in \ansin{H^1(\Omega)}$ is the \ansin{test function} for mechanical deformation and $f$ the linear functional corresponding to the force of external forces. The integral part with $v_M \boldsymbol{\Sigma}(u)\cdot n$ is often referred as a traction force which is applied on as a Neumann boundary condition on a surface. In our case, this traction force will be zero, as we neglect external gas pressure.
This is also everything we need for later implementation. We can use left hand side (lhs) and right hand side (rhs) function from FEniCS to extract the linear and bilinear forms. However for completion purpose we will conclude the linear and bilinear forms
\begin{align}
    a_T(T,v_T)&=\int_{\Omega}\rho c_p \dfrac{T_{n+1}-T_n}{\Delta t} \cdot v_T \ \dd x+\int_{\Omega}k\nabla T_{n+1} \cdot \nabla v_T \ \dd x+ \int_{\partial \Omega} hT_{n+1} v_T \ \dd A, \\
    a_M(u,v_m)&= \int_{\Omega} \boldsymbol{\Sigma}(u):\boldsymbol{\mathcal{E}}(v_M) \ \dd x,\\
    L_T(T)&= \int_{\partial \Omega} hT_\mathrm{e} v_T\dd A,\\
    L_M(u)&=\int_{\Omega} \boldsymbol{f} v_M \ \dd x.
\end{align}

The last thing we need to consider for the later implementation is the type of the finite element. We will use a second order Lagrange element for the displacement field and a first order Lagrange for the temperature. For details about finite element definition we refer to \citep{Saa21}, \citep{Rob12} and \citep{Log12}.

\subsection{Preprocessing and implementation}\label{subsec41}
We created our 3D data with the commercial engineering tool Autodesk Inventor. The full model can be seen in Figure \ref{fig:3d-figure} (a). Although the part is only an academic example, it is important that we place roundings at the otherwise sharp corners to get reasonably realistic results. To save memory and time we reduce the full model to a smaller one. Here we are using symmetry condition in x and y axis as well as using only one turbine blade as seen in Figure \ref{fig:3d-figure} (b). To use the  Cartesian coordinate system will simplify these boundary conditions afterwards. To mesh the 3D geometry, we must export the data to "stl" format, a common file format. Then we have to discretize the 3d model. We used "Gmsh" to do this. There we marked the whole volume model (tetrahedrons) as the computation domain $\Omega$ and necessary boundary surfaces (triangles) as the boundary conditions $\partial\Omega$. After the creation of the mesh we translated this to a mesh format we can use with FEniCS. We used "meshio" as a python module to convert the mesh into readable data structure for FEniCS. The overview in table formatting can be found in appendix.
\\
We implemented the thermoelasticty problem in FEniCS. This can be done by defining the problem as the weak formulation. We defined two function spaces for the two testfunctions from our weak form. We use two different element types in $\Omega$: the displacement field is approximated by Lagrange (FEniCS:Continous Galerkin) elements with order 2 and the temperature field is approximated by Lagrange elements of order one. In order to minimize computation time, we opted for the smallest possible element order that still yields reasonable stress results. The stress field has to be projected afterwards onto another space to relocate the stress to the nodes. Here we take a discontinuous Galerkin (DG) space of order 1 because the stress depends on the gradient of the displacement field which is of order 2. So we get at best a piecewise discontinuous results.
\\
The time is discretized in 20 equidistant steps. We first started with a logarithmic timestepping for better resolution in the first moments of the process. But it turns out that the time step should have some minimal value otherwise the solution of the heat equation shows some weird behavior. The values get negative despite the initial condition is nonnagtive as well as the boundary conditions are nonnegativ with positive time stepping with an backward Euler. This behavior is also seen in \citep{Cha15}. To preserve positivity we increased the time step. 
\\
To conclude our FEM, the input for this problem is the control of rotational speed and external temperature over 20 timesteps, which results in a 40 dimensional input. The model will give transient results of the stress and temperature over a time of 1800 seconds. We are especially interested in the maximum \ansin{von }Mises stress in each timestep.

\subsection{Definition of a starting guess control}
\begin{figure}[h]%
\centerline{\includegraphics[width=0.75\linewidth]{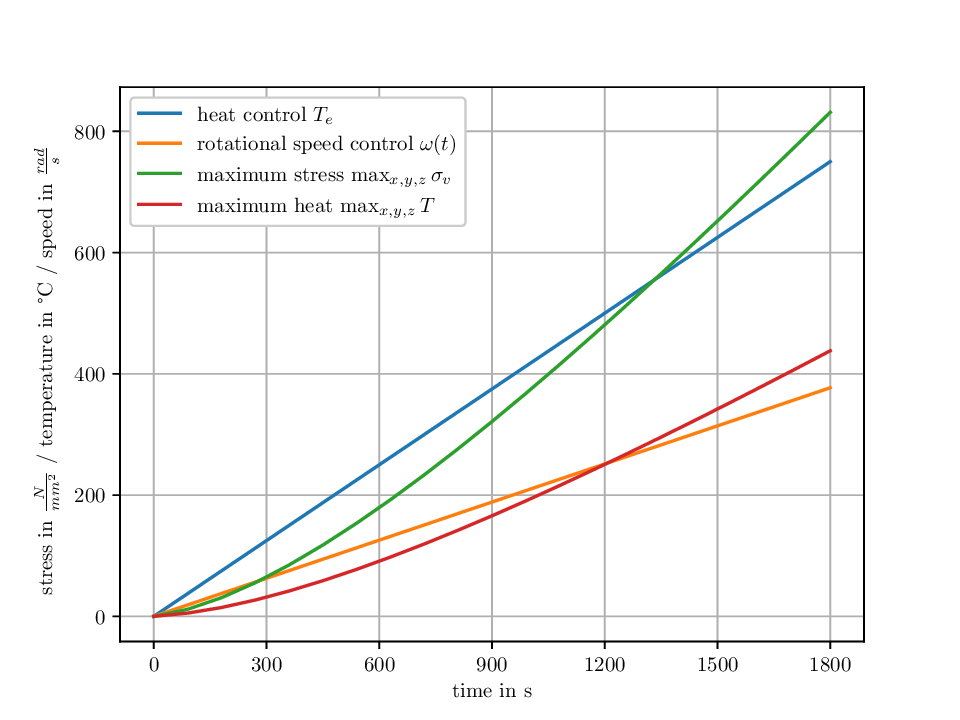}}
\caption{Initial guess of linear increasing speed and heat controls and solution of resulting maximum \ansin{von }Mises stress and maximum heat in the part over the whole simulation time. The linearly increasing control leads to the fulfillment of the constraints and seems to be a reasonable approach, since a slow heating naturally leads to a lower stress. But this initial start guess turns out to be non effective for gradient based optimization.}\label{fig:init_guess}
\end{figure}

The solution obtained from the finite element model corresponds to the input variables of external temperature and rotational speed. To establish a reasonable starting guess control, we adopt a linearly increasing rotation and temperature profile. This initial guess is illustrated in Figure \ref{fig:init_guess}. The outcome of a single finite element analysis entails determining the maximum \ansin{von }Mises stress (green) across the domain and the maximum heat (red) on the surface as shown in Figure \ref{fig:init_guess} with the two control variables heat control (blue) and rotational speed (orange). \ansin{This means we have a function (FEA) with 40 control variables, 20 external temperatures $T_e$ as indicated in blue and 20 rotational speeds $\omega$ as indicated in orange over the time domain. The response of the function is the maximum \ansin{von }Mises stress $\sigma_v$ and the maximum heat $T$ over the whole domain.}
The maximum \ansin{von }Mises stress occurs at the final time step and reaches a value exceeding $830 \frac{N}{mm^2}$. All imposed constraints and bounds are satisfied within this solution. This means, the maximum heat observed on the surface exceeds 400 °C, with the rotation rate reaching 60 Hz ($2\pi 60$ in plot) and the heat control reaching 750 °C at the conclusion of the analysis.

Nevertheless, in the subsequent chapter, we will demonstrate that this initial estimation is inadequate. We will present another initial guess, as indicated in Figure \ref{fig:iters}, that may not appear reasonable as a solid at first, but serves as a much better starting point for gradient-based optimization.

\subsection{Statistics and perfomance of the FEM model}\label{subsec43} 
In this finite element analysis (FEA), we used a mesh consisting of 1879 nodes, 1368 triangles, and 7339 tetrahedrons. The transient time step of the simulation was set to 90 seconds, and the simulation covered a total of 1800 seconds or 30 minutes with 20 time steps. When run on a cluster with hardware specs of 2x64 Core AMD EPYC Milan, 1024GB of RAM the same simulation took only about 1 minute. This demonstrates the power of utilizing a cluster for running FE simulations, as it can greatly reduce the computational time required and is even more relevant in the optimization. 

\section{Solving the optimal control problem}\label{sec5}
Solving optimal control problems involves finding the best set of control inputs to achieve a desired outcome within a system. One way to solve these problems is through the use of numerical optimization techniques, such as the Sequential Quadratic Programming (SQP) algorithm. In this case study, we will use an enhancement of the  SQP algorithm, the Sequential Least Squares Quadratic Programming (SLSQP) implemented in the Scipy library in Python which was developed in \citep{kra88}. We will take a look of the results of the optimization to understand how the control inputs affect the stress distribution within the part and understand the optimal solution that minimizes the overall stress.

In this chapter, we will explore the use of the SLSQP method for solving an optimal control problem. We will be using the Scipy library in Python, which is well suited for our problem. The implementation of Scipy is also able to deliver good results without a gradient which we will use and explain why we can not use automatic differentiation tools for the gradient.

\subsection{Optimal control with SLSQP method}\label{sec50}

\begin{figure}[h]%
\centerline{\includegraphics[width=0.95\linewidth]{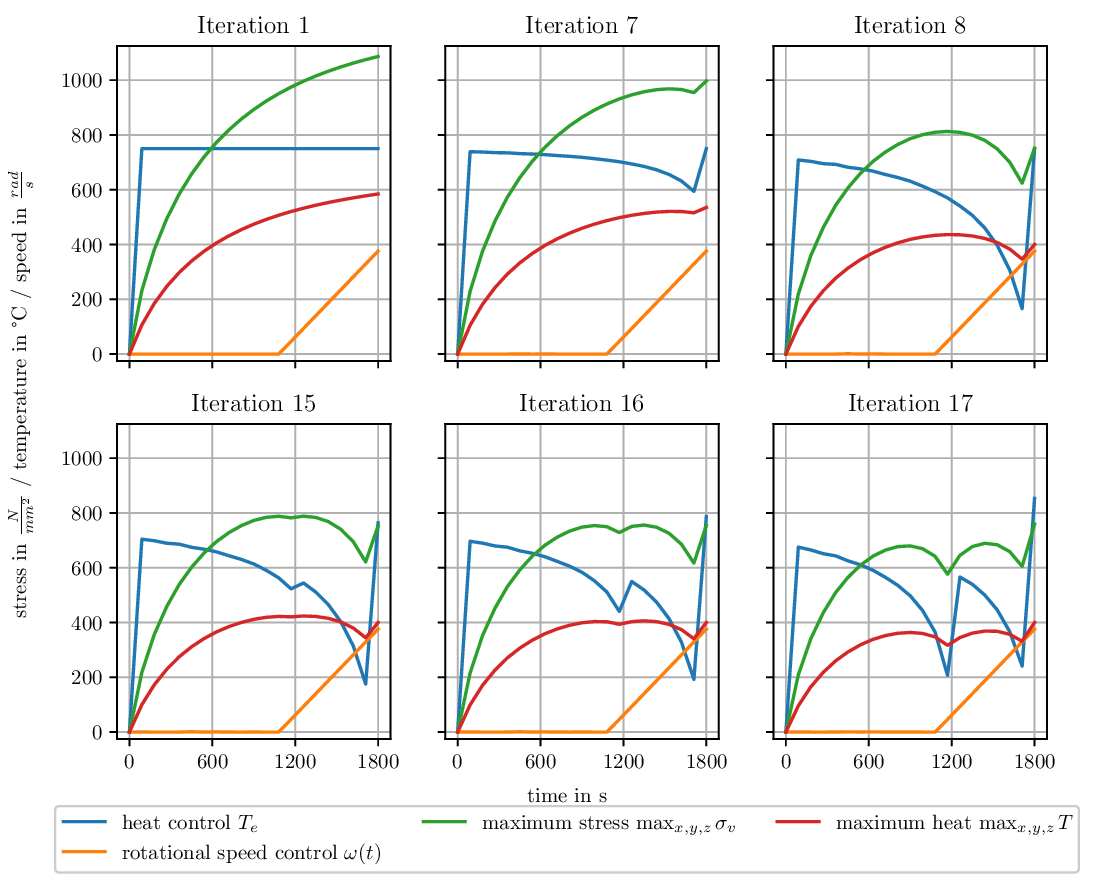}}
\caption{The maximum stress over the whole domain, as seen in green line, is decreasing from Iteration 1 to Iteration 17. In Iteration 1, another initial guess for the optimization procedure can be seen. A jump in heat control followed constant heat control and a speed control as late as possible will result in a better optimal control.}
\label{fig:iters}
\end{figure}

In this particular case, we will examine an optimal control problem that involves two control inputs: external heating and rotation. With a time step discretization of 20 and 2 input variables, we encounter a 40-dimensional optimization problem. The control variables are subject to specific limits: external heating must fall within the range of 0 °C to 1000 °C, while rotation must be between 0 and 60 Hz. However, it should be noted that the control temperature corresponds to the temperature of the working fluid. Furthermore, in addition to these bounds, the optimization algorithm includes additional constraints: the external heating should have a value of 750 °C at the final time step, and the rotation should be 60 Hz at the final time step. Additionally, we have incorporated a constraint to indicate the maximum temperature at the surface, which must reach 400 °C by the end of the simulation. 

The objective of this optimization is to minimize the maximum stress over all time steps and across the entire domain. Mathematically, this can be represented as shown in Equation \eqref{eq:opti}. The optimization algorithm iteratively adjusts the control inputs to minimize this objective, resulting in a decrease in the highest stress point, as illustrated by the green line in \ansin{Figure \ref{fig:iters}}.

Given that the SLSQP method is a gradient-based optimizer, providing a good initial guess for the control inputs is crucial. As mentioned earlier, we attempted using linearly increasing controls as the initial guess, but it turned out to be ineffective. A more effective initial solution involves maximum heat and minimal rotation, which helps in finding an optimal solution more quickly.

Various libraries offer automatic gradient calculation features. However, utilizing these automatic gradient procedures, such as Jax \citep{jax18} or Pyomo \citep{pyomo21}, requires the calculation to be implemented within the codebase of the respective libraries. Since we are using a C++-based code within FEniCS, compatibility is limited. Another option to obtain the gradient would be Dolfin-adjoint \citep{adjoint19}, similar to Tensorflow \citep{tf15}, which constructs a differentiable graph. However, this approach presents challenges as the objective function needs to be an integral over the entire domain, whereas our study specifically focuses on the maximum \ansin{von }Mises stress. Therefore, an integral over the entire domain is not useful in this research. So the gradient calculation is performed using finite differences in Scipy. We utilized the 2-point method, which is the default procedure when no gradient is provided to the optimization algorithm. The default step size is approximately $\sqrt{\textrm{machine precision}} \approx$ 1.49e-08, requiring one finite element method (FEM) call per optimization variable. Consequently, 40 FEM calls are made to obtain one gradient.

The SLSQP algorithm employs a BFGS approximation of the Hessian, as described by \citep{Powe78}, for constrained optimization problems. The optimal step size of the algorithm, analogous to Newton's method for nonlinear optimization problems, is $\alpha = 1$ near the local optimum. However, for vectors far from the optimum, the step size needs to be modified according to the penalty function described by \citep{Han77}, which includes a non-differentiable merit function. This function is substituted by a differentiable augmented Lagrange function presented by \citep{Schi80}.

\begin{figure}[h]%
\centerline{\includegraphics[width=0.75\linewidth]{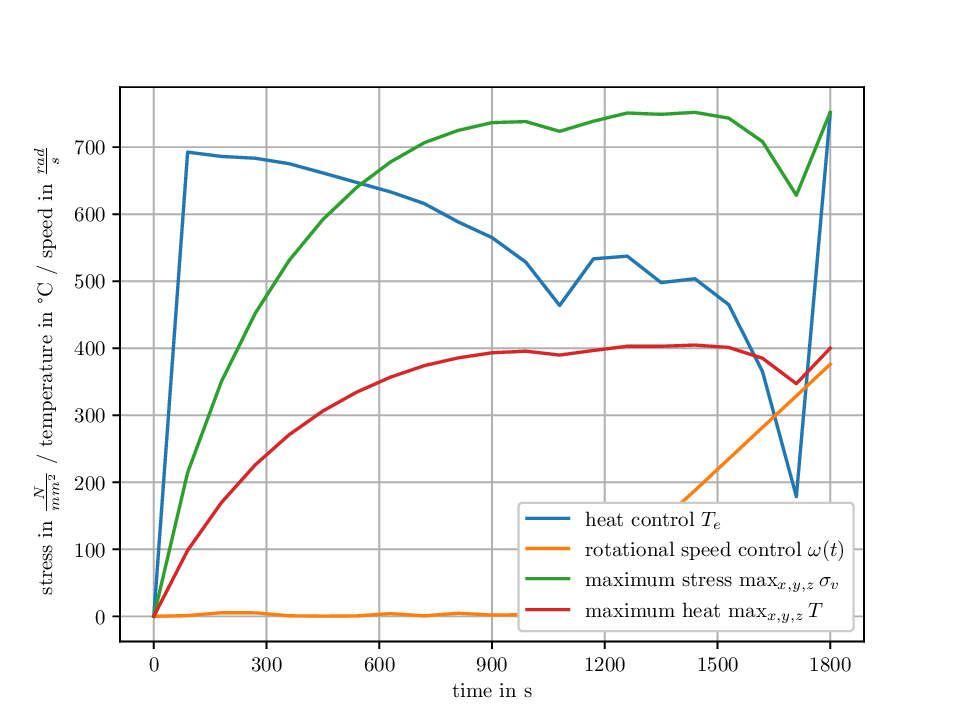}}
\caption{Optimal control results of Scipy's SLQSP in terms of speed and heat control. Heat control starts at a high value and then decreases as the stress increases. The stress increases and reaches a maximum, and the control tends to bring the stress as close to the maximum as possible, since the heat in the part can then increase as fast as possible. In the penultimate time step, both the \ansin{heat} control and the stress decrease to meet the boundary conditions in the final time steps while not exceeding a maximum stress threshold.}\label{fig:opt_con}
\end{figure}

In total, 179 iterations and over 7000 function evaluations (FEM simulation calls) were required to identify an optimal control with a function tolerance of $10^{-8}$. The optimal control is visualized in Figure \ref{fig:opt_con}. The SQP procedure invokes the FEM solver over 7000 times, resulting in a cumulative computation time of 7000 minutes or 120 hours. Therefore, the need for a fast FEM simulation is clearly evident. When the SQP method necessitates multiple restarts with different initial guesses, this time requirement becomes even more critical. 

\subsection{Physical Interpretation}\label{sec51}
In the present study, we observe a tendency in the optimal control strategy wherein the heating process is prioritized initially, followed by a delayed initiation of rotation. This sequential approach appears to be reasonable due to the nature of the influence exerted by rotational speed on \ansin{von }Mises stress, which is perceived as an offset. Conversely, the influence of temperature on stress diminishes once the component is fully heated. The optimal control tries for control to get a constant \ansin{von }Mises stress as the objective function is only influenced by the maximum over time. This results in an optimal use of time if the stress is nearly constant.

\begin{figure}[h]%
\centerline{\includegraphics[width=0.75\linewidth]{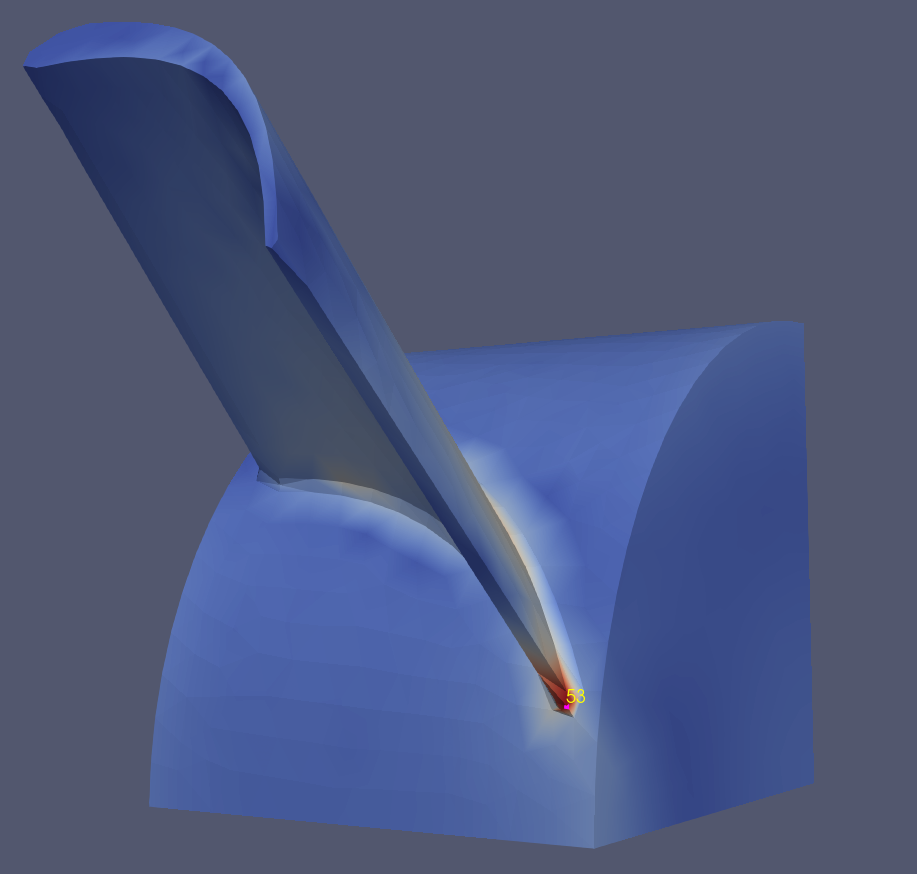}}
\caption{Position of the maximum \ansin{von }Mises stress in timestep 20. Only a small partition of the blade is under high mechanical stress. The highest stress is marked by a pink dot and the node number in yellow.}\label{fig:stress_pos}
\end{figure}

The position of the maximum stress is nearly at the same place. This is also one major advantage for gradient optimization when minimizing a maximum, because the the objective function will not take big jumps as response and remain continuous, which is not directly given by this objective function. The stress positions illustrated in Figure \ref{fig:stress_pos}. 

Moreover, our proposed solution leads to a reduction in maximum \ansin{von }Mises stress and increase of internal heat. Specifically, the \ansin{von }Mises stress is lowered from $830 \frac{N}{mm^2}$ \ansin{(see Figure \ref{fig:init_guess}, last time step)} to $750 \frac{N}{mm^2}$ \ansin{(see Figure \ref{fig:opt_con})}, resulting in a decrease of approximately 10\%. Additionally, we establish that minimizing rotation is desirable to mitigate the occurrence of unavoidable \ansin{von }Mises stress. \ansin{The break in \ansin{von }Mises stress and heat control in the penultimate time step is either explainable due to the evolution of the optimization (see Figure \ref{fig:iters}) where a "overshoot" can happen and never be reversed. This can occur due to the objective function which only states the maximum \ansin{von }Mises stress and therefore never change non-maxima. Another explanation is, that the constraints of this optimization problem at the last time step needs to be satisfied. The algorithm needs to reduce heat control in penultimate time step in order to satisfy constraints and not increase the maximum \ansin{von }Mises stress in the last time step.}

At first we did not bound the gradient of the rotational speed which resulted in a jump in the speed at the last timestep, which is technically not possible.

\section{Summary and Outlook}\label{sec6}
The  paper addresses the startup of a turbine through a transient finite element model incorporating thermoelasticity, with the heat equation, and rotational stress over 20 timesteps with two  controls that can be varied over time: the rotational speed and the external heat.  The aim is to minimize peaks in the \ansin{von }Mises stress through a control strategy of first heating and then rotating it as late as possible. This results in a convenient use of material and time in the sense that as much heat as possible is saved in the material before the rotational stress is applied. The slope of the rotational speed was limited by a constraint and a physically reasonable control was found. The \ansin{von }Mises stress peaks have been decreased by 10 \% which can lead to longer lifetime or safer components overall. \ansin{Up to the knowledge of the authors, for the first time, optimal control of an activation process with two control variables was successfully realized in a three-dimensional domain, leading to a 10\% reduction in peak \ansin{von }Mises stress.}

 \ansin{There are two types of limitation to this analysis. The first one concerns the abstracted physical behavior: we assumed the convective heat coefficient $h$ to be fixed and the inertia force is also neglected in this study which is not true in reality. This can be implemented in future work. The second limitation is} the need for faster simulation models when extending the optimization from one blade to a full turbine model. Creating a surrogate model can drastically decrease the computational time for one simulation. The full FEM model, or single FEM outputs, or the objective can be trained in such an approach. The surrogate model must be trained globally as well as in the optimal region. One can use latin hypercube sampling for global input parameters and an iterative approach with Bayesian optimization \citep{Sno12} to get new samples in the optimal region. Finally, also new methods such as physics-informed neural networks \citep{Rai19} could be used.

\section*{Disclosure statement}
The authors declare that no conflict of interest regarding financial or personal relationships exists that could have influenced this work.

\section*{Data availability statement}
A repository containing the optimal control script including the FEM solver and a mesh preprocessing can be found at: \url{https://github.com/nifri004/OC_Thermoelasticity}

\section*{Funding}
This paper was produced as part of the joint research project ML-Real-Time, which is funded by the German Federal Ministry of Education and Research (BMBF) through the VDI Technologiezentrum as project management agency of the BMBF for the funding program FHprofUnt 2018 (funding code. 13FH174PX8) based on a resolution of the German Bundestag and by Siemens Gas and Power GmbH \& Co. KG financially supported
H.G. also thanks the ’Forschungsvereinigung der Arbeitsgemeinschaft der Eisen und Metall verarbeitenden Industrie e.V.’ (AVIF No.A316) for their financial support.


\bibliographystyle{tfnlm}
\bibliography{opt_control_rota}

\newpage
\appendix
\section{}

\subsection{Workflow and used coefficents}
\begin{table}[h]
\centering
\begin{tabular}{lll}
\textbf{Workstep}  & \textbf{Used Program / Tool} & \textbf{File Format} \\
\hline
3D-CAD Data Generation           & Autodesk Inventor     & .stl       \\
\hline
Selecting Boundaries And calculation Area          & gmsh     & .geo       \\
\hline
General Meshing           & gmsh     & .mesh       \\
\hline
Convert Mesh For FEniCS         & meshio (python)     & .xmdf, .h5      \\
\hline
Run FEM Simulation         & FEniCS (python)     &        \\
\hline
Run Optimization          & Scipys SQP (python)     &       \\
\hline
Show 3D Results         & Paraview     &   .xmdf, .h5   \\
\end{tabular}
\caption{Workflow for optimizing a control a 3d model in FEM. Listed are the required work steps, the resulting file formats and the used programs and tools}
\label{tbl:workflow}
\end{table}

\begin{table}[h]
\centering
\begin{tabular}{lll}
\textbf{Name}  & \textbf{Abbreviation} & \textbf{Value and Unit} \\[8pt]
\hline
E-Modulus           & $E$       &  210e9 $\frac{N}{m^2}$ \\ [8pt]
\hline
lame coefficient 1         & $\lambda$       &  $\frac{\nu}{1-2\nu}\cdot\frac{E}{1+\nu}$ \\ [8pt]
\hline
lame coefficient 2          & $\mu$       &  $\frac{1}{2}\cdot\frac{E}{1+\nu}$ \\ [8pt]
\hline
Poisson number          & $\nu$       &  0.3 \\ [8pt]
\hline
density         & $\rho$       &  $8050 \frac{kg}{m^3}$ \\ [8pt]
\hline
Thermal expansion coefficient        & $\alpha$       &  13.5e-6 $\frac{1}{K}$ \\ [8pt]
\hline
Specific heat        & $c_p$       &  $ 420 \frac{J}{kG K}$  \\ [8pt]
\hline
Thermal conductivity       & $k$       &  $36 \frac{W}{m K}$ \\ [8pt]
\hline
convective heat coefficient       & $h$       &  $20 \frac{W}{m^2 K}$ \\ [8pt]
\hline
\end{tabular}
\caption{Used constant and coefficients}
\label{tbl:coefs}
\end{table}

\newpage
\subsection{Additional information about stress location}

\begin{figure}[h]
\centering
\subfigure[Timestep 4 with maximal \ansin{von }Mises stress marked with dot and node number.]{%
\resizebox*{6.5cm}{!}{\includegraphics{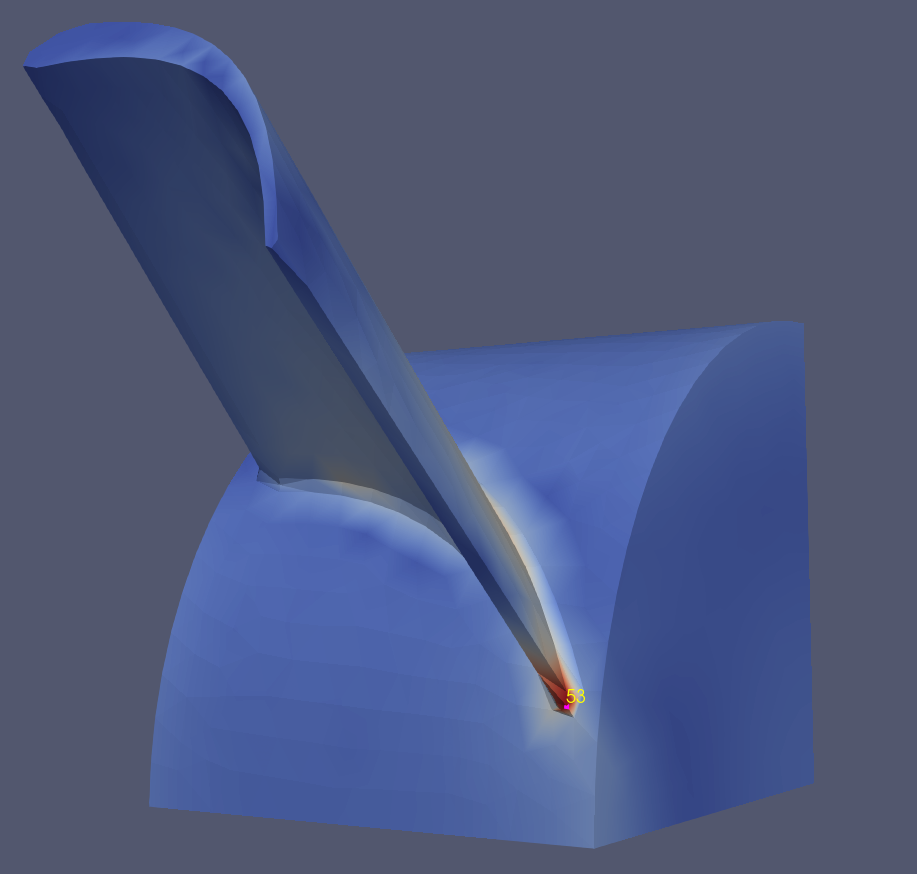}}}\hspace{5pt}
\subfigure[Timestep 5 with maximal \ansin{von }Mises stress marked with dot and node number.]{
\resizebox*{6cm}{!}{\includegraphics{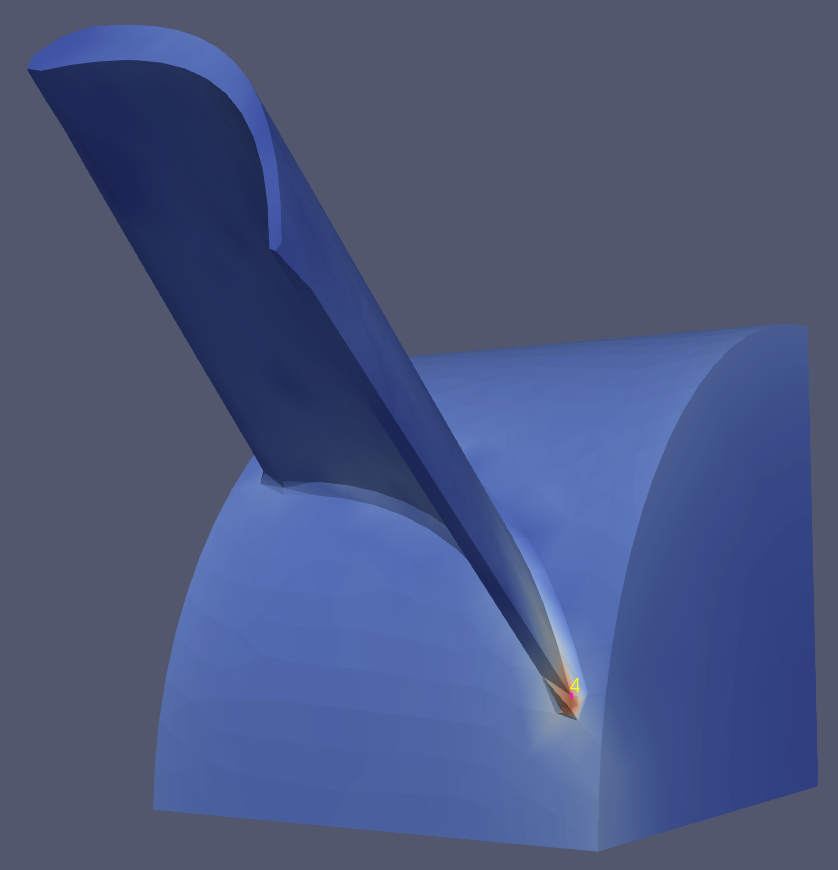}}}
\caption{Von Mises stress results for timestep 4 and 5. Location of maximal stress is moving between this points. This is important for our gradient based optimizer, otherwise it could result in jumps in stress.} 
\label{fig:app:3d-figure}
\end{figure}

\begin{table}[h]
\centering
\begin{tabular}{lllllllllll}
\textbf{timestep} & 1 & 2 &3 &4 &5 &6 &7 &8 &9 &10 \\
\hline
\textbf{node number} & 53 & 53 & 53 & 53 & 4 & 4 & 4 & 4 & 4 & 4\\
\hline \hline
\textbf{timestep} &11 &12 &13 &14 &15 &16 &17 &18 & 19 & 20 \\
\hline
\textbf{node number} & 4 & 4 & 4 & 4 & 4 & 4 & 4 & 53 & 53 & 53       \\
\end{tabular}
\caption{Location of maximal van Mises stress}
\label{tbl:max_stress_pos}
\end{table}

\end{document}